\documentclass[a4paper, final, notitlepage, 10pt]{article}
\usepackage{latexsym,bm}
\usepackage{amsthm}
\usepackage{amsfonts}
\usepackage{amsmath}
\usepackage[all]{xy}
\usepackage[final]{graphics}
\usepackage{fancyhdr}
\usepackage{verbatim}
\usepackage{geometry}
\usepackage{xcolor}
\usepackage{enumerate}

\newtheorem{thm}{Theorem}[section]
\newtheorem{cor}[thm]{Corollary}
\newtheorem{lem}[thm]{Lemma}
\newtheorem{prop}[thm]{Proposition}

\theoremstyle{definition}

\newtheorem{rem}[thm]{Remark}
\newtheorem{exa}[thm]{Example}

\numberwithin{equation}{section}

\newcommand{\thmref}[1]{Theorem~\ref{#1}}

\newcommand{\lemref}[1]{Lemma~\ref{#1}}
\newcommand{\propref}[1]{Proposition~\ref{#1}}
\newcommand{\corref}[1]{Corollary~\ref{#1}}
\newcommand{\exaref}[1]{Example~\ref{#1}}
\newcommand{\remref}[1]{Remark~\ref{#1}}

\renewcommand{\O}{\mathcal{O}}

\renewcommand{\L}{\mathcal{L}}

\newcommand{\PP}{\mathbb{P}}

\newcommand{\bpi}{\overline{\pi}}
\newcommand{\e}{\varepsilon}
\newcommand{\nequiv}{\stackrel{num}{\sim}}
\newcommand{\tr}{\mathrm{tr}}

\renewcommand{\i}{\mathrm{i}}

\title{Notes on automorphisms of surfaces of general type \\ with $p_g=0$ and $K^2=7$}
\author{Yifan Chen}
\date{}
\begin{document}
\maketitle
\renewcommand{\thefootnote}{\fnsymbol{footnote}}
\footnotetext{Yifan Chen: Academy of Mathematics and Systems Science, Chinese Academy of Sciences, No.~55 Zhongguancun East Road, Haidian District, Beijing, 100190, P.~R.~China}
\footnotetext{email address: chenyifan1984@gmail.com}

\begin{abstract}
Let $S$ be a smooth minimal complex surface of general type with $p_g=0$ and $K^2=7$.
We prove that any involution on $S$ is in the center of the automorphism group  of $S$.
As an application, we show that the automorphism group of an Inoue surface with $K^2=7$
is isomorphic to $\mathbb{Z}_2^2$ or $\mathbb{Z}_2 \times \mathbb{Z}_4$.
We construct a $2$-dimensional family of Inoue surfaces with automorphism groups isomorphic to $\mathbb{Z}_2 \times \mathbb{Z}_4$.
\end{abstract}

\section{Introduction}
The birational automorphism groups of projective varieties are extensively studied.
Nowadays we know that,
for a projective variety  of general type $X$ over an algebraically closed field of characteristic zero,
the number of birational automorphisms of $X$ is bounded by $c_d \cdot \mathrm{vol}(X, K_X)$,
where $c_d$ is a constant which only depends on the dimension $d$ of $X$,
and $\mathrm{vol}(X, K_X)$ is the volume of the canonical divisor $K_X$ (cf.~\cite{HMX}).
Furthermore, we know that $c_1=42$ and $c_2=42^2$ from the classical Hurwitz theorem and Xiao's theorem
(cf.~\cite{xiao1} and \cite{xiao2}).
However, even in low dimensions,
it is usually nontrivial to calculate the automorphism groups of explicit varieties of general type
(for example, see \cite{fakeprojectiveplane1}, \cite{fakeprojectiveplane2}, \cite{fakeprojectiveplane3} and \cite{godeauxauto}).

We focus on automorphisms of minimal smooth complex surfaces of general type with $p_g=0$ and $K^2=7$.
Involutions on such surfaces have been studied in \cite{leeshin} and  \cite{rito}.
All the possibilities of the quotient surfaces and the fixed loci of the involutions are listed.
In order to find new examples, we have tried to classify such surfaces with commuting involutions in \cite{commuting}
and succeeded in constructing a new family of surfaces in \cite{chennew}.
We briefly recall the main results in \cite{commuting}.
Throughout the article,
$S$ denotes a minimal smooth surface of general type with $p_g=0$ and $K^2=7$ over $\mathbb{C}$.
\begin{thm}(\cite[Theorem~1.1, Theorem~2.9 and Section~4]{commuting})\label{thm:classification}
Assume that the automorphism group $\mathrm{Aut}(S)$ contains a subgroup $G=\{1, g_1, g_2, g_3\}$,
which is isomorphic to $\mathbb{Z}_2^2$.
Let $R_{g_i}$ be the divisorial part of the fixed locus of the involution $g_i$ for $i=1, 2, 3$.
Then the canonical divisor $K_S$ is ample and $R_{g_i}^2=-1$ for $i=1, 2, 3$.
Moreover, there are only three numerical possibilities for the intersection numbers
$(K_SR_{g_1}, K_SR_{g_2}, K_SR_{g_3})$:
$(a)\ (7, 5, 5)$, $(b)\ (5, 5, 3)$ and $(c)\ (5, 3, 1)$.
The intersection numbers $(R_{g_1}R_{g_2}, R_{g_1}R_{g_3}, R_{g_2}R_{g_3})$ have the following values:
$(a)\ (5, 9, 7)$, $(b)\ (7, 5, 1)$ and $(c)\ (1, 3, 1)$, respectively.
\end{thm}
In the above theorem, we adopt the convention that $K_SR_{g_1} \ge K_SR_{g_2} \ge K_SR_{g_3}$,
$R_{g_1}R_{g_2} \le R_{g_1}R_{g_3}$ in case (a) and $R_{g_1}R_{g_3} \ge R_{g_2}R_{g_3}$ in case (b).
Actually,
we have completely classified the surfaces in case (a) and case (b) in \cite{commuting}.
But we do not know any example of the surfaces in case (c).
One may ask whether there are noncommutative involutions on $S$.
Here we give a negative answer.
\begin{thm}\label{thm:center}
If $\alpha$ is an involution of $S$, then $\alpha$ is contained in the center of $\mathrm{Aut}(S)$.
\end{thm}
We prove the above theorem in Section~3.
The key step is \thmref{thm:commutative} which shows that any two involutions on $S$ commute.
\thmref{thm:commutative} also has the following corollary.
\begin{cor}\label{cor:allinvolutions}
Assume that $(S, G)$ is a pair satisfying the assumption of \thmref{thm:classification}.
Then $g_1, g_2$ and $g_3$ are exactly all the involutions of $\mathrm{Aut}(S)$.
\end{cor}
The corollary immediately implies that
if $\mathrm{Aut}(S)$ contains a nontrivial subgroup which is isomorphic to $\mathbb{Z}_2^r$,
then $r=1$ or $r=2$.
We remark that there are surfaces of general type with $p_g=0$, $K^2=8$ and
their automorphism groups contain subgroups which are isomorphic to $\mathbb{Z}_2^3$
(cf.~\cite[Example~4.2--4.4]{bicanonical1}).

As an application,
we calculate the automorphism groups of the surfaces in the case (a)  of \thmref{thm:classification}.
These surfaces are those constructed by M.~Inoue in \cite{inoue}
and they are the first examples of surfaces of general type with $p_g=0$ and $K^2=7$.
They can be described as finite Galois $\mathbb{Z}_2^2$-covers of the $4$-nodal cubic surface
(see~\exaref{exa:inoue}, which is from \cite[Example~4.1]{bicanonical1}).
\begin{thm}\label{thm:inoueauto}
Let $S$ be an Inoue surface.
Then $\mathrm{Aut}(S) \cong \mathbb{Z}_2^2$ or $\mathrm{Aut}(S) \cong \mathbb{Z}_2 \times \mathbb{Z}_4$.
If $S$ is a general Inoue surface, then $\mathrm{Aut}(S) \cong \mathbb{Z}_2^2$.
\end{thm}

Inoue surfaces form a $4$-dimensional irreducible connected component in the Gieseker moduli space of
canonical models of surfaces of general type (cf.~\cite{inouemfd}).
The proof of \thmref{thm:inoueauto} actually shows that
$\mathrm{Aut}(S) \cong \mathbb{Z}_2^2$ for $S$ outside a $2$-dimensional irreducible closed subset of this connected component (see \remref{rem:2dim}).
We also exhibit a $2$-dimensional family of Inoue surfaces with
$\mathrm{Aut}(S) \cong \mathbb{Z}_2 \times \mathbb{Z}_4$ (see Section~5).
They are finite Galois $\mathbb{Z}_2 \times \mathbb{Z}_4$-covers of a $5$-nodal weak Del Pezzo surface of degree two,
which is the minimal resolution of one node of the $6$-nodal Del Pezzo surface of degree two.

\section{Preliminaries}
\subsection{Fixed point formulae}
Let  $X$ be a smooth projective surface over the complex number field.
We only consider surfaces with $p_g(X)=q(X)=0$.
In this case, $X$ has Picard number $\rho(X)=10-K_X^2$ by Noether's formula and Hodge decomposition.
Also the expotential cohomology sequence gives $\mathrm{Pic}(X) \cong H^2(X, \mathbb{Z})$.
Poincar\'e duality implies that the intersection form on
$\mathrm{Num}(X):=\mathrm{Pic}(X)/\mathrm{Pic}(X)_{\mathrm{Tors}}$ is unimodular.

Assume that $X$ has a nontrivial automorphism $\tau$.
Denote by $\mathrm{Fix}(\tau)$ the fixed locus of $\tau$.
Let $k_{\tau}$ be the number of isolated fixed points of $\tau$
and let $R_\tau$ be the divisorial part of $\mathrm{Fix}(\tau)$.
Then $R_\tau$ is a disjoint union of irreducible smooth curves.
We denote by $\tau^* \colon H^2(X, \mathbb{C}) \rightarrow  H^2(X, \mathbb{C})$
the induced linear map on the second singular cohomology group (note that $H^k(X, \mathbb{C})\\=0$ for $k=1, 3$).
The following proposition follows directly from the Topological and Holomorphic Lefschetz Fixed Point Formulae
(cf.~\cite{AS}, Page 567; see also \cite[Lemma~4.2]{manynodes}).
The automorphism $\tau$ is called an involution if it is of order $2$.

\begin{prop}\label{prop:lefschetz}
If $\tau$ is an involution,
then
$k_\tau=K_XR_\tau+4\ \text{and}\ \tr(\tau^*)=2-R_\tau^2.$
If $\tau$ is of order $3$,
then
$k_\tau=r_1+r_2=\tr(\tau^*)+2+K_XR_\tau+R_\tau^2\ \text{and}\ r_1+2r_2=6+\frac 32K_XR_\tau-\frac {R_\tau^2}{2},$
where $r_j$ is the number of isolated fixed points of $\tau$ of type $\frac 13(1, j)$ for $j=1, 2$.
\end{prop}

Throughout this article,
we denote by $S$ a smooth minimal complex surface of general type with $p_g=0$ and $K_S^2=7$.
Then $\rho(S)=3$ and $S$ contains at most one $(-2)$-curve
(this follows from Poincar\'e duality; cf.\cite[Lemma~2.5]{commuting}).
Here an $m$-curve (for $m \le 0$) on a smooth surface
stands for an irreducible smooth rational curve with self intersection number $m$.
\begin{lem}\label{lem:oneinvolution}(See also the table in \cite{leeshin})
Let $\tau$ be an involution on $S$.
Then $K_SR_\tau \in \{1, 3, 5, 7\}$ and $R_\tau^2= \pm 1$.
If $R_\tau^2=1$, then $K_S$ is ample and $R_\tau$ is irreducible with $K_SR_\tau=3$.
\end{lem}
\proof
For $R_\tau^2= \pm 1$, see the proof of \cite[Proposition~3.6]{inouebloch}.
According to \cite[Lemma~3.2 and Proposition~3.3~(v)]{involution},
$k_\tau$ is an odd integer and $k_\tau \le 11$.
So $K_SR_\tau \in \{1, 3, 5, 7\}$ by \propref{prop:lefschetz}.

Assume that $R_\tau^2=1$.
If $K_S$ is not ample, then $S$ has a unique $(-2)$-curve $C$.
The intersection number matrix of $K_S, R_\tau$ and $C$ has determinant
$-14+2(K_SR_\tau)^2-7(R_\tau C)^2$.
The determinant equals $0$,
for otherwise, the Chern classes of $K_S, R_\tau$ and $C$ form a basis of $H^2(S, \mathbb{C})$
and they are $\tau^*$-invariant, a contradiction to $\tr(\tau^*)=2-R_\tau^2=1$ by \propref{prop:lefschetz}.
It follows that $K_SR_\tau=7$ and $(R_\tau C)^2=12$.
This is impossible. So $K_S$ is ample.

The algebraic index theorem gives $(K_SR_{\tau})^2 \ge K_S^2 R_{\tau}^2=7$ and thus $K_SR_{\tau} \in \{3, 5, 7\}$.
Let $\pi_\tau \colon S \rightarrow \Sigma_\tau:=S/<\tau>$ be the quotient morphism.
We have $K_S=\pi_\tau^*(K_{\Sigma_\tau})+R_\tau$.

If $K_SR_{\tau}=5$, then
$k_{\tau}=K_SR_{\tau}+4=9$ and $K_{\Sigma_{\tau}}^2=\frac 12 (K_S-R_{\tau})^2=-1$.
So $\Sigma$ has $9$ nodes.
If $\Sigma_\tau$ has Kodaira dimension $\kappa(\Sigma_{\tau}) \ge 0$,
then the minimal resolution $W_{\tau}$ of $\Sigma_\tau$ has Picard number $11$
and it contains $9$ disjoint $(-2)$-curves.
By \cite[Proposition~4.1]{manynodes},
$W_\tau$ is minimal.
This contradicts $K_{W_{\tau}}^2=-1$.
So $\kappa(\Sigma_{\tau})=-\infty$
and $W_{\tau}$ is a rational surface.
This contradicts \cite[Theorem~3.3]{manynodes}.
Hence $K_SR_{\tau} \not=5$.

In the same manner we see that $K_SR_{\tau} \not=7$ (see also \cite{bicanonical2}).
So $K_SR_{\tau}=3$.
Because $K_S$ is ample and $R_{\tau}$ is a disjoint union of smooth irreducible curves,
the algebraic index theorem shows that $R_{\tau}$ is irreducible.
\qed

\subsection{Abelian covers}
We briefly recall some facts from the theory of abelian covers from \cite{pardini}.
Assume that $\pi \colon X \rightarrow Y$ is a finite abelian cover between projective varieties
with $X$ normal and $Y$ smooth.
Let $\mathfrak{S}$ be the Galois group of $\pi$
and let $\mathfrak{S}^*$ be the group of characters of $\mathfrak{S}$.
Then the action of $\mathfrak{S}$ induces a splitting:
$\pi_\ast(\O_X)=\oplus_{\chi \in \mathfrak{S}^*} \L_{\chi}^{-1}$,
where $\L_{\chi} \in \mathrm{Pic}(Y)$ and $\L_1=\O_X$.
For each nontrivial cyclic subgroup $\mathfrak{C}$ of $\mathfrak{S}$ and each generator $\psi \in \mathfrak{C}^*$,
there is a unique effective divisor $D_{\mathfrak{C}, \psi}$ of $Y$ associated to the pair $(\mathfrak{C}, \psi)$.
The cover $\pi$ is determined by $\L_{\chi}$ and $D_{\mathfrak{C}, \psi}$
with some specified relations (cf.~\cite[Theorem~2.1]{pardini}).
We mainly apply this theory when $\mathfrak{S} \cong \mathbb{Z}_2 \times \mathbb{Z}_4$ or
$\mathfrak{S} \cong \mathbb{Z}_2^2$.

We set up some notation and conventions.
Denote by $H=<g_1>\times<g>$ a group isomorphic to $\mathbb{Z}_2 \times \mathbb{Z}_4$,
where $g_1, g$ are generators of $H$, $g_1$ is of order $2$ and $g$ is of order $4$.
Denote by $H^*=<\chi> \times <\rho>$ the group of characters of $H$,
where $\chi(g_1)=-1, \rho(g)=\i$ and $\chi(g)=\rho(g_1)=1$.
The group $H$ contains a unique subgroup $G=\{1, g_1, g_2, g_3\}$ which is isomorphic $\mathbb{Z}_2^2$,
where $g_2=g^2$ and $g_3=g_1g_2$.
Denote by $\chi_i \in G^*$ the nontrivial character orthogonal to $g_i$ for $i=1, 2, 3$.

When $\mathfrak{S}=G$,
we simply set $\L_i:=\L_{\chi_i}$ and $\Delta_i:=D_{<g_i>, \psi}$,
where $\psi$ is the unique nontrivial character of $<g_i>$.
Similarly, when $\mathfrak{S}=H$, we set $D_i:=D_{<g_i>, \psi}$ for $1 \not =\psi \in <g_i>^*$.
For the cyclic group $<g> \cong \mathbb{Z}_4$,
we set $D_{g, \pm \i}:=D_{<g>, \psi}$ for $\psi \in <g>^*$ with $\psi(g)=\pm \i$.
We adopt similar convention for the cyclic group $<g_1g> \cong \mathbb{Z}_4$.

In what follows,
the indices $i \in \{1, 2, 3\}$ should be understood as residue classes modulo $3$.
Also linear equivalence and numerical equivalence between divisors are denoted by $\equiv$ and $\nequiv$, respectively.

\begin{prop}(cf.~\cite{singularbidouble}, \cite[Theorem~2.1 and Corollary~3.1]{pardini})\label{prop:abeliancover}
Let $\pi \colon X \rightarrow Y$ be a finite abelian cover between projective varieties.
Assume that $X$ is normal and $Y$ is smooth.
\begin{enumerate}[\upshape (a)]
    \item If the Galois group of $\pi$ is  $G$, then $\pi$ is determined by the following data:
          \begin{align}
          2\L_i \equiv \Delta_{i+1}+\Delta_{i+2},\ \L_i+\Delta_i \equiv \L_{i+1}+\L_{i+2}\ \text{for}\ i=1, 2, 3, \label{eq:G-cover}
          \end{align}
          where $\L_i, \Delta_i$ are divisors of $Y$, $\Delta_i$ is effective
          and $\Delta:=\Delta_1+\Delta_2+\Delta_3$ is reduced.
    \item If the Galois group of $\pi$ is $H$, then $\pi$ is determined by the following reduced data
          (see \cite[Proposition~2.1]{pardini}):
          \begin{align}
          2\L_{\chi} &\equiv D_1+D_3+D_{g_1g, \i}+D_{g_1g, -\i},\label{eq:Hcover}\\
          4\L_{\rho} &\equiv 2D_2+2D_3+D_{g, \i}+3D_{g, -\i}+D_{g_1g,\i}+3D_{g_1g,-\i}, \nonumber
          \end{align}
          where $\L_{\chi}, \L_{\rho}$ and $D_1, \ldots, D_{g_1g,-\i}$ are divisors of $Y$,
          $D_1, \ldots, D_{g_1g,-\i}$ are effective and
          \begin{align*}
          D:=D_{1}+D_{2}+D_{3}+D_{g,\i}+D_{g,-\i}+D_{g_1g, \i}+D_{g_1g, -\i}\ \text{is reduced.}
          \end{align*}
\end{enumerate}
\end{prop}

\section{Two involutions commute}
We first deduce \corref{cor:allinvolutions} and \thmref{thm:center}
from the following theorem.
\begin{thm}\label{thm:commutative}
Let $S$ be a smooth minimal surface of general type with $p_g(S)=0$ and $K_S^2=7$.
Assume that $\mathrm{Aut}(S)$ contains two distinct involutions $\alpha$ and $\beta$.
Then $\alpha\beta=\beta\alpha$.
\end{thm}

\proof[Proof of \corref{cor:allinvolutions}]
On the contrary, suppose that $\alpha$ is an involution of $\mathrm{Aut}(S)$ other than $g_1, g_2, g_3$.
\thmref{thm:commutative} implies $<\alpha, g_1, g_2> \cong \mathbb{Z}_2^3$.
This group contains seven subgroups of order $4$, say $G_1, \ldots, G_7$.
Each pair $(S, G_j)$ must satisfy one of the three numerical possibilities in \thmref{thm:classification}.
However, this is impossible because any two of these seven subgroups have a common involution.
Hence $g_1, g_2$ and $g_3$ are exactly all the involutions of $\mathrm{Aut}(S)$.
\qed

\proof[Proof of \thmref{thm:center}]
We may assume that $\mathrm{Aut}(S)$ contains at least two involutions.
These two involutions generate a subgroup  $G \cong \mathbb{Z}_2^2$ by \thmref{thm:commutative}.
We still denote by $g_1, g_2$ and $g_3$ the involutions of $G$.
Let $\tau$ be any automorphism of $S$.
\corref{cor:allinvolutions} gives $\tau G \tau^{-1}=G$.
Since $\tau(R_{g_i})=R_{\tau g_i \tau^{-1}}$,
we have $K_SR_{g_i}=K_SR_{\tau g_i\tau^{-1}}$ and
$R_{g_i}R_{g_{i+1}}=R_{\tau g_i\tau^{-1}}R_{\tau g_{i+1}\tau^{-1}}$ for $i=1, 2, 3$.
From this observation and \thmref{thm:classification},
we conclude that $\tau g_i\tau^{-1}=g_i$ for $i=1, 2, 3$ and complete the proof.
\qed\smallskip

The remaining of this section is devoted to prove \thmref{thm:commutative}.
{\textbf{We assume by contradiction that $\mathbf{\alpha\beta \not= \beta\alpha}$.}}
We will deduce a contradiction through a sequence of lemmas and propositions.
We use the same notation as Section~2.
Recall that $\tr(\alpha^*)=2-R_{\alpha}^2$, $R_\alpha^2=\pm 1$
and $k_\alpha=K_SR_\alpha+4$
(see \propref{prop:lefschetz} and \lemref{lem:oneinvolution}).
\begin{lem}\label{lem:oddorder}
The order of $\alpha \beta$ is an odd integer.
\end{lem}
\proof
    Assume by contradiction that the order of $\alpha\beta$ is $2k$ and $k \ge 2$.
    Let $\gamma :=(\alpha \beta)^k = (\beta \alpha )^k$.
    Then $\gamma$ is an involution and
    $\gamma \alpha = (\alpha \beta)^k \alpha = \alpha (\beta \alpha)^k = \alpha \gamma$.
    Therefore $<\gamma, \alpha> \cong \mathbb{Z}_2^2$.
    Then $R_{\alpha}^2=R_{\gamma}^2=-1$ by \thmref{thm:classification}.
    Similarly, $\gamma \beta = \beta \gamma$ and $R_{\beta}^2=R_{\gamma}^2=-1$.
    So $\tr(\alpha^*)=\tr(\beta^*)=3$.

    Let $\iota:=\alpha \beta \alpha$.
    Note that $\alpha, \beta$ and $\iota$ are three distinct involutions in $\mathrm{Aut}(S)$ and
    \begin{align}
    \alpha(R_{\iota})=R_\beta,\ \alpha(R_{\beta})=R_{\iota},\ \alpha(R_{\alpha})=R_\alpha \label{eq:conjugate1}
    \end{align}
    Recall that $\dim H^2(S, \mathbb{C})=\rho(S)=3$.
    Now $c_1(R_{\alpha}), c_1(R_{\beta})$ and $c_1(R_{\iota})$ are not a basis of $H^2(S, \mathbb{C})$,
    for otherwise, \eqref{eq:conjugate1} implies $\tr(\alpha^*)=1$,
    which is a contradiction to $\tr(\alpha^*)=3$.

    So the intersection number matrix of $R_{\alpha}, R_{\beta}$ and $R_{\iota}$ has determinant zero.
    That is $2x^2y+2x^2+y^2-1=0$,
    where $x:=R_{\alpha}R_{\iota}=R_{\alpha}R_{\beta}$ (see \eqref{eq:conjugate1}) and $y:=R_{\beta}R_{\iota}$.
    It follows that $x=0, y=1$ and the nontrivial linear relation among
    $c_1(R_\alpha), c_1(R_{\beta})$ and $c_1(R_\iota)$ is $c_1(R_{\beta})+c_1(R_{\iota}) = 0$.
    This contradicts the fact that the divisor $R_{\beta}+R_\iota$ is strictly effective.
    Hence the order of $\alpha\beta$ is an odd integer. \qed\smallskip

Recall that our aim is to deduce a contradiction from the assumption $\alpha\beta \not= \beta\alpha$.
According to the previous lemma, from now on,
{\textbf{we may assume that the order $\mathbf{r}$ of $\mathbf{\alpha\beta}$ is an odd prime.}}
In fact, if $r=p(2t+1)$ for some prime $p \ge 3$ and some integer $t > 0$,
then $\alpha':=(\alpha\beta)^t\alpha$ and $\beta':=(\beta\alpha)^t\beta$
are involutions and the order of $\alpha'\beta'$ is $p$.
We may replace $\alpha, \beta$ by $\alpha', \beta'$ and continue our discussion.

The subgroup $<\alpha, \beta>$ of $\mathrm{Aut}(S)$ is isomorphic to the dihedral group of order $2r$.
Let $D_r$ denote this subgroup.
Since $r$ is a prime, all the involutions in $D_r$ are pairwise conjugate and
$D_r$ has exactly one nontrivial normal subgroup $<\alpha\beta>$, which is the commutator subgroup.
Any irreducible linear representation of $D_r$ has dimension at most two,
and any irreducible $2$-dimensional representation of $D_r$ is isomorphic to the matrix representation given by
$$\left(\begin{array}{cc}
    0 & 1\\
    1 & 0
    \end{array}\right)
    \\\text{and}\
    \left(\begin{array}{cc}
    c & 0\\
    0 & c^{-1}
    \end{array}\right)
    \\\text{for some}\ c \not= 1\ \text{and}\ c^r=1.$$
\begin{lem}\label{lem:nontrivialaction}With the same assumption as above, we have
\begin{enumerate}[\upshape (a)]
    \item the canonical $K_S$ is ample;
    \item the curves $R_\alpha$ and $R_\beta$ generate a pencil $|F|$ of curves with $F^2=1$ and $K_SF=3$,
          and $|F|$ has a simple base point $p$;
    \item the group $D_r$ acts faithfully on $|F|$.
\end{enumerate}
\end{lem}
\proof
    Assume that the order of $\alpha\beta$ is $r=2k+1$ for $k \ge 1$.
    Set $\gamma:=\alpha(\beta \alpha)^k = \beta (\alpha \beta)^k$.
    Then $\alpha, \beta$ and $\gamma$ are three distinct involutions and they are pairwise conjugate.
    Therefore $K_SR_{\alpha}=K_SR_{\beta}=K_SR_{\gamma}$ and $R_{\alpha}^2=R_{\beta}^2=R_{\gamma}^2$.
    Since $\gamma \alpha = \beta \gamma$ and $\gamma \beta = \alpha \gamma$,
    \begin{align}
    \gamma (R_{\alpha})=R_{\beta},\ \gamma (R_{\beta})=R_{\alpha},\ \gamma (R_{\gamma})=R_{\gamma} \label{eq:conjugate2}
    \end{align}

    We claim that $R_{\gamma}^2=R_{\alpha}^2=R_{\beta}^2=1$.
    Otherwise, as in the proof of \lemref{lem:oddorder},
    we could deduce a contradiction by calculating the determinant of the intersection number matrix of
    $R_{\alpha}, R_{\beta}$ and $ R_{\gamma}$ and by calculating $\tr(\gamma^*)$.

    Then (a) follows from \lemref{lem:oneinvolution}.
    \lemref{lem:oneinvolution} also gives $K_SR_\alpha=K_SR_\beta=3$.
    The the algebraic index theorem  implies
    $(R_{\alpha}+R_{\beta})^2 \le \frac {6^2}{7}$ and thus $R_{\alpha}R_{\beta} \le 1$.
    Since $R_{\alpha}^2=R_{\beta}^2=1$,
    the equality holds and $R_{\alpha} \nequiv R_{\beta}$.
    Similarly, we have $R_{\gamma} \nequiv R_{\alpha}.$

    Let $p$ be the unique intersection point of $R_{\alpha}$ and $R_{\gamma}$.
    Then \eqref{eq:conjugate2} implies that $R_{\alpha}, R_{\beta}$ and $R_{\gamma}$
    pairwise intersect transversely at the point $p$.
    Recall that $\mathrm{Pic}(S) \cong H^2(S, \mathbb{Z})$ and
    $\mathrm{Num}(S)=\mathrm{Pic}(S)/\mathrm{Pic}(S)_{\mathrm{Tors}}$.
    Let $m$ be the smallest positive integer such that $mR_{\alpha} \equiv mR_{\gamma} \equiv mR_{\beta}$.
    Let $\varepsilon \colon \tilde{S} \rightarrow S$ be the blowup at $p$,
    let $E$ be the exceptional curve
    and let $\tilde{R_{\alpha}}$ be the strict transform of $R_\alpha$, etc.
    Then $|m\tilde{R_{\alpha}}|$ induces a fibration $f \colon \tilde{S} \rightarrow \mathbb{P}^1$
    and $m\tilde{R_{\gamma}}, m\tilde{R_{\alpha}}$ and $m\tilde{R_{\beta}}$ are fibers of $f$.

    The fibration $f$ has $E$ as a $m$-section.
    If $m \ge 2$, we easily obtain a contradiction by
    applying the Hurwitz formula for $f|_E \colon E \rightarrow \mathbb{P}^1$.
    Therefore $m=1$, $R_{\alpha} \equiv R_{\gamma} \equiv R_{\beta}$ and
    $h^0(S, \O_S(R_{\alpha}))=2$.
    And (b) is proved.

    For (c), first note that $p$ is a fixed point of $D_r$.
    So $D_r$ acts faithfully on the tangent space $\mathrm{T}_pS$ of $S$ to the point $p$.
    According to the discussion before the lemma, this action is irreducible
    and the corresponding action of $D_r$ on $\mathbb{P}(\mathrm{T}_pS)$ is faithful.
    Since $F^2=1$, $p$ is a smooth point of $F$
    and thus $\mathrm{T}_pF$ is a $1$-dimensional linear subspace of $\mathrm{T}_pS$ for any $F \in |F|$.
    From this, we conclude that $D_r$ acts faithfully on $|F|$.\qed\smallskip

Because $D_r$ acts faithfully on $|F| \cong \mathbb{P}^1$,
every automorphism has exactly two invariant curves in $|F|$.
For every involution $\gamma \in D_r$,
one of the two $\gamma$-invariant curves in $|F|$ is $R_{\gamma}$.
Denote the other one by $F_{\gamma}$.
Then $F_{\gamma}$ contains the seven isolated fixed points of $\gamma$.
Denote by $F_0$ one of the two $\alpha\beta$-invariant curves in $|F|$.
Then the other one is $\alpha(F_0) (=\beta(F_0))$ and $\mathrm{Fix}(\alpha\beta) \subseteq F_0 \cup \alpha(F_0)$.
We shall show that $F_0$ is not $2$-connected.
But first we need the following lemma about the action of $D_r$ on the singular cohomology group.

\begin{lem}\label{lem:quotientpicard}
The automorphism $\alpha\beta$ acts trivially on $H^2(S, \mathbb{C})$.
In particular, the quotient surface $S/D_r$ has Picard number $2$.
\end{lem}
\proof
We have seen that
$\alpha, \beta$ and thus $D_r$ act trivially on the $2$-dimensional linear subspace generated $c_1(K_S)$ and $c_1(F)$.
Because $H^2(S, \mathbb{C})$ is $3$-dimensional and $\alpha\beta$ is contained in the kernel of
any $1$-dimensional representation of $D_r$,
$\alpha\beta$ acts trivially on $H^2(S, \mathbb{C})$.
Hence the invariant subspace of $H^2(S, \mathbb{C})$ for the $D_r$-action is $2$-dimensional and
$S/D_r$ has Picard number $2$.
\qed\smallskip

We analysis the members of the pencil $|F|$, which are not $2$-connected.
This will help us to determine the base locus of the linear system $|K_S+F|$ in the proof of \propref{prop:orderthree}
and to find a basis of $\mathrm{Num}(S)$.
We continue to use the fact that $S$ has Picard number $3$.

\begin{lem}\label{lem:non2connected}
Assume that $|F|$ contains a curve which is not $2$-connected.
Then
\begin{enumerate}[\upshape (a)]
    \item the curves in $|F|$ which are not $2$-connected are exactly the $\alpha\beta$-invariant curves $F_0$ and $\alpha(F_0)$;
    \item $F_0=A+B$, where $A$ and $B$ are irreducible curves, and
          $K_SA=2, K_SB=1, A^2=0, B^2=-1$ and $AB=1$. Moreover, $A$ contains the base point $p$ of $|F|$.
\end{enumerate}
\end{lem}

\proof
Assume that $A+B \in |F|$, $A>0$, $B>0$ and $AB \le 1$.
Because $K_S$ is ample and $K_SF=3$,
we may assume $K_SA=2$ and $K_SB=1$.
Then $B$ is irreducible.
The algebraic index theorem implies $A^2 \le 0$ and $B^2 \le -1$.
In particular, by \lemref{lem:quotientpicard}, $\alpha\beta(B).B=B^2<0$
and thus $\alpha\beta(B)=B$.
Hence $A+B$ is one of the $\alpha\beta$-invariant curves $F_0$ and $\alpha(F_0)$.

Because $A^2+B^2=F^2-2AB \ge -1$,
the argument above yields $A^2=0$, $B^2=-1$ and $AB=1$.
Then $FA=1$ and $FB=0$.
So the simple base point $p$ of $|F|$ belongs to $A$.
It remains to show that $A$ is irreducible.
Assume by contradiction that $A$ is reducible.
Because $K_SA=2$ and $K_S$ is ample,
$A=A_1+A_2$, $K_SA_1=K_SA_2=1$ and both $A_1$ and $A_2$ are irreducible.
We may  assume $p \in A_1$ and $p \not \in A_2$.
Then $FA_1=1, FA_2=0$ and $A_2^2<0$.
The adjunction formula gives $A_2^2=-1$ or $A_2^2=-3$.
The intersection number matrix of $K_S, F$ and $A_2$ has determinant $-2A_2^2-1$.
Since the intersection form on $\mathrm{Num}(S)$ is unimodular,
we get $A_2^2=-1$.
Since $F^2=1$ and $F(A_2+B)=0$,
the algebraic index theorem implies $A_2B=0$.
But then $S$ contains four disjoint curves $B, \alpha(B), A_2$ and $\alpha(A_2)$,
all with self intersection number $(-1)$,
a contradiction to $\rho(S)=3$.\qed\smallskip

The following proposition determines the order of the automorphism $\alpha\beta$.
\begin{prop}\label{prop:orderthree}
The automorphism $\alpha\beta$ is of order $3$.
Moreover, $F_0$ is not $2$-connected,
where $F_0$ is as in \lemref{lem:non2connected}.
\end{prop}

\proof
Let $F$ be  any curve in $|F|$.
The long exact sequence of cohomology groups associated to the  exact sequence
$0 \rightarrow \O_S(K_S) \rightarrow \O_S(K_S+F) \rightarrow \O_{F}(K_{F}) \rightarrow 0$
shows that $h^0(S, \O_S(K_S+F))=h^0(F, \O_{F}(K_{F}))=p_a(F)=3$
and the trace of $|K_S+F|$ on $F$ is complete.
Thus $|K_S+F|$ defines a rational map $h \colon S \dashrightarrow \mathbb{P}^2$
and $h$ is defined on $F$ whenever $|K_F|$ is base point free.
In particular,
$h$ is defined on the smooth curve $R_\alpha\ (\in |F|)$ and $h(R_\alpha)$ is the canonical image of $R_\alpha$.
The same statement holds by replacing $\alpha$ by $\beta$.

Because there is a $D_r$-linearization on $\O_S(K_S+F)$, the rational map $h$ is $D_r$-equivalent.
Therefore $h(R_\alpha)$ is contained in the fixed locus of the action of $\alpha$ on $\mathbb{P}^2$.
Note that an involution on $\mathbb{P}^2$ has a line and a point as the fixed locus.
It follows that $\alpha$ acts trivially on $\mathbb{P}^2$ because $h(R_\alpha)$ is a conic curve or a quartic curve.
Similarly, $\beta$ and thus $D_r$ act trivially on $\mathbb{P}^2$.
Therefore $h \colon S \dashrightarrow \mathbb{P}^2$ factors through the quotient morphism $S \rightarrow S/D_r$.

Note that $K_S$ is ample, $F$ is nef and $(K_S+F)^2=14$.
First assume that $h$ is a morphism.
Then it is finite and it has degree $14$.
We thus get $|D_r|=\deg h$ and $r=7$.
It follows that the induced morphism $h' \colon S/D_r \rightarrow \mathbb{P}^2$ is an isomorphism.
So the invariant linear subspace of $H^2(S, \mathbb{C})$ for the $D_r$-action is isomorphic to
$H^2(\mathbb{P}^2, \mathbb{C})$, which is $1$-dimensional.
This contradicts \lemref{lem:quotientpicard}
and thus $h$ is not a morphism.

We now analysis the base locus of $h$.
If $F$ is $2$-connected, $|K_{F}|$ is base point free by \cite[Theorem~3.3]{embedding}
and $h$ is defined on $F$.
Hence the base locus of $|K_S+F|$ is contained in the curves of $|F|$, which are not $2$-connected.
According to \lemref{lem:non2connected},
$|F|$ contains exactly two such curves $F_0=A+B$ and $\alpha(F_0)$.
Similar arguments as above show that the trace of $|K_S+A|$ (respectively $|K_S+B|$) on
$A$ (respectively $B$) is complete.
Since $p_a(A)=2$ and $p_a(B)=1$,
$|K_A|$ and $|K_B|$ are base point free by \cite[Theorem~3.3]{embedding}.
Because $|K_S+F| \supseteq |K_S+A|+B, |K_S+B|+A$,
$|K_S+F|$ has exactly two base points $q:=A \cap B$ and $\alpha(q)=\alpha(A) \cap \alpha(B)$.

Therefore $h$ is a finite morphism outside the base locus and $\deg h=(K_S+F)^2-2=12$.
Since $h$ factors through $S/D_r$, we have $|D_r|=6$ and $r=3$.
\qed\smallskip

It is easy to check that $F, A$ and $\alpha(A)$ generate $\mathrm{Num}(S)$
and $K_S \nequiv F+A+\alpha(A)$.
We shall show that $K_S$ is indeed linearly equivalent to $F+A+\alpha(A)$,
and deduce a contradiction to $p_g(S)=0$ and complete the proof of \thmref{thm:commutative}.
For this purpose, we turn to the quotient surface $S/D_r$
and analysis $\mathrm{Fix}(\alpha\beta)$.

\begin{prop}\label{prop:fixedlocus}
The automorphism $\alpha\beta$ has $B \cup \alpha(B)$ as the divisorial part of the fixed locus
and it has five isolated fixed points $p, q_1, q_2, \alpha(q_1)$ and $\alpha(q_2)$,
where $q_1$ and $q_2$ are contained in $A$.
Each isolated fixed point of $\alpha\beta$ is of type $\frac {1}{3} (1, 2)$.
\end{prop}
\proof
We have seen that $F_0$ and $\alpha(F_0)$ are $\alpha\beta$-invariant
and $\mathrm{Fix}(\alpha\beta) \subseteq F_0 \cup \alpha(F_0)$.
Moreover, the curves $A, \alpha(A), B$ and $\alpha(B)$ are $\alpha\beta$-invariant.
Also note that a point $q$ is a fixed point (respectively an isolated fixed point) of $\alpha\beta$
if and only if so is the point $\alpha(q)(=\beta(q))$.

We claim that neither $A$ nor $\alpha(A)$ is contained in $\mathrm{Fix}(\alpha\beta)$.
Otherwise, both $A$ and $\alpha(A)$ are contained in $\mathrm{Fix}(\alpha\beta)$.
Since $A \cap \alpha(A)=p$,
this contradicts the fact that the divisorial part of $\alpha\beta$ is a disjoint union of smooth curves.
The claim is proved.

Now assume by contradiction that $B$ is not contained in $\mathrm{Fix}(\alpha\beta)$.
Then nor is $\alpha(B)$ and $\mathrm{Fix}(\alpha\beta)$ consists of isolated fixed points.
Then $\mathrm{Fix}(\alpha\beta)$ has five fixed points by \propref{prop:lefschetz} and \lemref{lem:quotientpicard}.
Three of these points are $p, q:=A \cap B$ and $\alpha(q)$.
Denote the other two by $p_1(\in F_0)$ and by $\alpha(p_1)$.
We must have $p_1 \in B$.
Otherwise, the nontrivial automorphism $\alpha\beta|_B$ has exactly one fixed point $q$,
which is a smooth point of $B$ since $AB=1$.
This is impossible because $p_a(B)=1$.
Therefore $p_1 \in B$.
It follows that $\alpha\beta|_A$ has exactly two fixed points $p$ and $q$,
which are smooth points of $A$.
Note that $A$ has at most two singular points since $p_a(A)=2$.
Because the singular locus of $A$ is $\alpha\beta$-invariant and $\alpha\beta|_A$ has order $3$,
we conclude that $A$ is indeed smooth.
However, the Hurwitz formula shows that $\alpha\beta|_A$ has either one or four fixed points,
a contradiction.

So $B$ and $\alpha(B)$ are contained in $\mathrm{Fix}(\alpha\beta)$.
In particular, $B$ and $\alpha(B)$ are smooth curves.
Then $\mathrm{Fix}(\alpha\beta)\setminus \{B \cup \alpha(B)\}$ consists of five isolated fixed points
and each fixed point is of type $\frac 13(1, 2)$ by \propref{prop:lefschetz} and \lemref{lem:quotientpicard}.
These points must be contained in $A \cup \alpha(A)$.
\qed\smallskip

Now we are able to describe the quotient map $\pi \colon S \rightarrow Y:=S/D_r$,
where $D_r=\{1, \alpha, \beta, \gamma, \alpha\beta, \beta\alpha\}$ and $\gamma:=\alpha\beta\alpha=\beta\alpha\beta$.
The divisorial parts and isolated fixed points of cyclic subgroups of $D_r$ are as follows
(see the discussion before \lemref{lem:quotientpicard} and \propref{prop:fixedlocus}):
\begin{align*}
 \text{cyclic subgroups}   && \text{divisorial part}  && \text{isolated fixed points}\\
 <\alpha> (\text{resp.}\ <\beta>, <\gamma>) \cong\mathbb{Z}_2 && R_{\alpha}\ (\text{resp.}\ R_{\beta}, R_{\gamma})             && \text{$7$ points on}\ F_{\alpha}\ (\text{resp.}\ F_{\beta}, F_{\gamma})\\
 <\alpha\beta>\cong\mathbb{Z}_3 &&B, \alpha(B)  &&p, q_1, q_2, \alpha(q_1), \alpha(q_2)
\end{align*}

Note that $p$ is the unique point with the stabilizer $D_r$.
From the action of $D_r$ on the tangent space $\mathrm{T}_pS$ (see the proof of \lemref{lem:nontrivialaction}~(c)),
it is easily seen that $\pi(p)$ is a smooth point of $Y$.
We conclude that $Y$ has seven nodes and two $A_2$-singularities $\pi(q_1)$ and $\pi(q_2)$.
In particular, $Y$ is Gorenstein.
The ramification formula gives
\begin{align}
K_S=\pi^*K_Y+R_{\alpha}+R_{\beta}+R_{\gamma}+2B+2\alpha(B)
   \equiv \pi^*K_Y+3F+2B+2\alpha(B) \label{eq:KSKQ}
\end{align}
and thus $K_Y^2=\frac 16 (K_S-3F-2B-2\alpha(B))^2=-3$.

Let $B'=\pi(B)$.
Then $B'$ is contained in the smooth locus of $Y$.
Note that $B'$ is a smooth elliptic curve and $\pi^*B'=3B+3\alpha(B)$.
So $B'^2=-3$ and $K_QB'=3$.
Since $(-K_QB')=K_Y^2B'^2$,
$-K_Y \nequiv B'$ by \lemref{lem:quotientpicard}.
This implies that $H^0(mK_Y)=0$ for $m \ge 1$.
As the quotient of $S$, $Y$ has irregularity $q(Y)=0$.
Therefore $Y$ is a rational surface.
Note that linear equivalence and numerical equivalence between divisors are the same on a smooth rational surface.
Since $Y$ contains only rational double points and $B'$ is contained in the smooth locus of $Y$,
we have $-K_Y\equiv B'$ indeed.
Then by \eqref{eq:KSKQ},
\begin{align*}
K_S \equiv \pi^*(-B')+3F+2B+2\alpha(B) \equiv (-3B-3\alpha(B))+3F+2B+2\alpha(B) \equiv F+A+\alpha(A).
\end{align*}
We obtain a contradiction to $p_g(S)=0$ and complete the proof of \thmref{thm:commutative}.

\section{Inoue Surfaces}
As mentioned in the introduction,
Inoue surfaces are the first examples of surfaces of general type with $p_g=0$ and $K^2=7$ (cf.~\cite{inoue}).
Here we describe them as finite Galois $\mathbb{Z}_2^2$-covers of the $4$-nodal cubic surface,
following \cite[Example~4.1]{bicanonical1}.
At the end of this section, we prove \thmref{thm:inoueauto}.
\begin{exa}\label{exa:inoue}
Let $\sigma \colon W \rightarrow \mathbb{P}^2$ be the blowup of the six vertices $p_1,p_2,p_3,p_1',p_2',p_3'$
of a complete quadrilateral on $\mathbb{P}^2$.
\begin{figure}[h]
\begin{centering}
\includegraphics{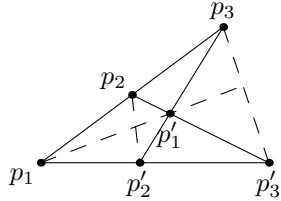}
\caption{\footnotesize{Configurations of the points $p_1, \ldots, p_3'$}\label{figure1}}
\end{centering}
\end{figure}
Denote by $E_i$ (respectively $E_i'$) the exceptional curve of $W$ over $p_i$ (respectively $p_i'$)
and denote by $L$ the pullback of a general line by $\sigma$.
Then $\mathrm{Pic}(W)=\mathbb{Z}L\oplus \oplus_{i=1}^3(\mathbb{Z}E_i\oplus \mathbb{Z}E_i')$.

The surface $W$ has four disjoint $(-2)$-curves. They are the proper transforms of the four sides of
the quadrilateral and their divisor classes are
\begin{align*}
Z_i \equiv L-E_i-E_{i+1}'-E_{i+2}',\  Z \equiv L-E_1-E_2-E_3.
\end{align*}
Let $\eta \colon W \rightarrow \Sigma$ be the morphism contracting there curves.
Then $\Sigma$ is the $4$-nodal cubic surface.

Let $\Gamma_1,\Gamma_2$ and $\Gamma_3$ be the proper transforms of the three diagonals of the quadrilateral,
i.e., $\Gamma_i \equiv L-E_i-E_i'$ for $i=1, 2, 3$.
Note that they are exactly the $(-1)$-curves which are disjoint from any $(-2)$-curve.
For each $i=1,2,3,$
$W$ has a pencil of rational curves $|F_i|:=|2L-E_{i+1}-E_{i+2}-E_{i+1}'-E_{i+2}'|$.
Observe that $-K_W \equiv \Gamma_1+\Gamma_2+\Gamma_3 \equiv \Gamma_i+F_i$ for $i=1,2,3$.

We define three effective divisors on $W$
\begin{align}
\Delta_1:=\Gamma_1+F_2+Z_1+Z_3,\ \Delta_2:=\Gamma_2+F_3,\ \Delta_3:=\Gamma_3+F_1+F_1'+Z_2+Z  \label{eq:inouedata}
\end{align}
We require that $F_i$ ($i=1, 2, 3$) and $F_1'$ are smooth $0$-curves such that the divisor
$\Delta: =\Delta_1+\Delta_2+\Delta_3$ has only nodes.
It is directly to show that there are divisors $\L_1, \L_2$ and $\L_3$
satisfying \eqref{eq:G-cover} in \propref{prop:abeliancover}.
Then there is a smooth finite $G$-cover $\bpi \colon V \rightarrow W$
branched on the divisors $\Delta_1, \Delta_2$ and $ \Delta_3$.
The (set theoretic) inverse image of a $(-2)$-curve under $\bpi$ is a disjoint union of two $(-1)$-curves.
Let $\e \colon V  \rightarrow S$ be the blowdown of these eight $(-1)$-curves.
Then there is a finite $G$-cover $\pi \colon S \rightarrow \Sigma$ such that
the following diagram \eqref{diag} commutes.
\begin{align}\label{diag}
\xymatrix{
    V \ar^{\e}"1,2" \ar_{\bpi}"2,1" & S \ar^{\pi}"2,2" \\
    W \ar^{\eta}"2,2" & \Sigma
}
\end{align}
The surface $S$ is a smooth minimal surface of general type with $p_g(S)=0$ and $K_S^2=7$.
It is called an Inoue surface.
When the curves $F_1, F_1', F_2$ and $F_3$ vary,
we obtain a $4$-dimensional family of Inoue surfaces.
\end{exa}

\begin{lem}\label{lem:identityinoue}
Let $W$ be as in \exaref{exa:inoue}.
\begin{enumerate}[\upshape (a)]
    \item Let $\alpha$ be an automorphism on $W$.
          If the induced map $\alpha^* \colon H^2(W, \mathbb{C}) \rightarrow H^2(W, \mathbb{C})$ is the identity,
          then $\alpha=\mathrm{Id}_W$.
    \item Let $\alpha_{\mathbb{P}^2}$ be the involution on $\mathbb{P}^2$ such that
          $\alpha_{\mathbb{P}^2}(p_k)=p_k'$ for $k=1, 3$.
           It induces an involution $\alpha_0$ on $W$.
           Then $\mathrm{Fix}(\alpha_0)$ consists of the $(-1)$-curve $\Gamma_2$
           and three isolated fixed points $\Gamma_1 \cap \Gamma_3$, $E_2 \cap F_3^*$ and $E_2' \cap F_3^*$,
           where $F_3^*$ is the
           unique smooth $\alpha_0$-invariant curve in the pencil $|F_3|=|2L-E_1-E_1'-E_2-E_2'|$.
\end{enumerate}
\end{lem}
\proof
For (a), the assumption implies that the $(-1)$-curves $E_i$ and $E_i'$ ($i=1, 2, 3$) are $\alpha$-invariant.
So $\alpha$ comes from an automorphism on $\mathbb{P}^2$ which has $p_1, \ldots, p_3'$ as fixed points
and thus it is the identity morphism.

For (b), note that $\mathrm{Fix}(\alpha_{\mathbb{P}^2})=\overline{p_2p_2'} \cup \{p_{13}:=\overline{p_1p_1'} \cap \overline{p_3p_3'}\}$ and $\sigma(\mathrm{Fix}(\alpha_0))=\mathrm{Fix}(\alpha_{\mathbb{P}^2})$.
Because $\sigma^{-1}(\overline{p_2p_2'})=E_2 \cup E_2' \cup \Gamma_2$ and the divisorial part of $\alpha_0$ is smooth,
$\mathrm{Fix}(\alpha_0)$ has $\Gamma_2$ as the divisorial part.
Then $\alpha_0$ has three isolated fixed points by \propref{prop:lefschetz}.
The point $\sigma^{-1}(p_{13})=\Gamma_1 \cap \Gamma_3$ is an isolated fixed point of $\alpha_0$.
Note that $\alpha_0$ induces a nontrivial action on the pencil $|F_3|$.
So $|F_3|$ contains exactly two $\alpha_0$-invariant curves $\Gamma_1+\Gamma_2$ and $F_3^*$.
Since $E_2$ and $E_2'$ are also $\alpha_0$-invariant, the intersection points $E_2 \cap F_3^*$ and $E_2' \cap F_3^*$
are isolated fixed points of $\alpha_0$.
\qed

\proof[Proof of \thmref{thm:inoueauto}]

Let $\tau \in \mathrm{Aut}(S)$.
By \thmref{thm:center}, $\mathrm{Fix}(g_i)$ is $\tau$-invariant for $i=1, 2, 3$ and
$\tau$ induces an automorphism $\alpha_\Sigma$ on the quotient surface $\Sigma = S/G$.
So the branch locus $\pi(\mathrm{Fix}(g_i))$ (for $i=1, 2, 3$) of $\pi \colon S \rightarrow \Sigma$
is $\alpha_\Sigma$-invariant.

Assume that $\tau \not \in G$, i.e., $\alpha_\Sigma \not =\mathrm{Id}_\Sigma$.
The automorphism $\alpha_\Sigma$ lifts to the minimal resolution $W$ of $\Sigma$.
Denote by $\alpha$ the induced automorphism on $W$.
Then $\Delta_1, \Delta_2, \Delta_3$ (see the diagram\eqref{diag}) are $\alpha$-invariant
because $\Delta_i$ is the inverse image of $\pi(\mathrm{Fix}(g_i))$
under the morphism $\eta \colon W \rightarrow \Sigma$.
These divisors are given by \eqref{eq:inouedata}.
It follows that the $(-1)$-curves $\Gamma_1, \Gamma_2, \Gamma_3$, the $0$-curves $F_2, F_3$
and the curves $F_1+F_1', Z_1+Z_3, Z_2+Z$ are $\alpha$-invariant.
Note that the Chern classes of $\Gamma_1, \Gamma_2, \Gamma_3$ and the Chern classes of $Z_1, Z_2, Z_3, Z$
generate $H^2(S, \mathbb{C})$.
The argument above implies that $(\alpha^2)^*=(\alpha^*)^2$ is the identity morphism.
Then $\alpha$ is an involution by \lemref{lem:identityinoue}.

Since $F_i\equiv \Gamma_{i+1}+\Gamma_{i+2}$,
the fibration $f_i \colon W \rightarrow \mathbb{P}^1$ induced by $|F_i|$ is $\alpha$-equivalent for $i=1, 2, 3$.
Note that $f_2$ has three singular fibers $\Gamma_1+\Gamma_3, Z_1+2E_2'+Z_3$ and $Z_2+2E_2+Z$.
According to the discussion above, these three fibers are $\alpha$-invariant.
Because any nontrivial automorphism on $\mathbb{P}^1$ has at most two fixed points,
$\alpha$ respects the fibration $f_2$, i.e., $f_2=f_2\alpha$.
In particular, $E_2$ and $E_2'$ are $\alpha$-invariant.

Note that $f_1$ has three singular fibers $\Gamma_2+\Gamma_3, Z_1+2E_1+Z$ and $Z_2+2E_1'+Z_3$.
If $f_1=f_1\alpha$, then all the $(-2)$-curves $Z_1, Z_2, Z_3$ and $Z$ are $\alpha$-invariant
since $\alpha$ also respects $f_2$.
Then $\alpha^*$ is the identity morphism and so is $\alpha$ by \lemref{lem:identityinoue},
a contradiction to our assumption.
So $\alpha$ induces a nontrivial action on $|F_1| \cong \mathbb{P}^1$.
Since the singular $\Gamma_2+\Gamma_3$ is $\alpha$-invariant,
$\alpha$ must permute the other two singular fibers of $f_1$.
Hence $\alpha(E_1)=E_1'$ and $\alpha(E_1')=E_1$.
Similarly, by considering the action of $\alpha$ on $|F_3|$, we see that
$\alpha(E_3)=E_3'$ and $\alpha(E_3')=E_3$.

We conclude that $\alpha$ is the involution $\alpha_0$ in \lemref{lem:identityinoue}.
We actually prove that if $\mathrm{Aut}(S) \not =G$, then $\mathrm{Aut}(S)/G \cong <\alpha_0>$,
and in the Equation~\eqref{eq:inouedata},
the curve $F_3$ in $\Delta_2$ is indeed the curve $F_3^*$ in \lemref{lem:identityinoue}~(b)
and $F_1'=\alpha_0(F_1)$ in $\Delta_3$.
Combining with \thmref{thm:center} and \corref{cor:allinvolutions},
we complete the proof of \thmref{thm:inoueauto}
\qed

\begin{rem}\label{rem:2dim}
When the curves $F_1$ and $F_2$ vary, the Inoue surfaces corresponding to the branch divisors \eqref{eq:inouedata} with $F_3=F_3^*$ and $F_1'=\alpha_0(F_1)$ form a $2$-dimensional irreducible closed subset
of the total $4$-dimensional family of Inoue surfaces.
Also \lemref{lem:identityinoue} shows that $W/<\alpha_0>$ has three nodes.
Moreover,
it contains three $(-2)$-curves in the smooth locus and these curves are the images of $Z_1+Z_3$, $Z_2+Z$ and $\Gamma_2$
under the quotient map from $W$ to $W/<\alpha_0>$.
This observation motivates us to construct some special Inoue surfaces in the next section.
\end{rem}

\section{Special Inoue surfaces}
We construct a $2$-dimensional family of Inoue surfaces with automorphism groups
isomorphic to $\mathbb{Z}_2 \times \mathbb{Z}_4$.
We use the notation in Subsection~2.2.

Let $q, q_1,q_2,q_3,q_1',q_2'$ and $q_3'$ be seven points on $\PP^2$ with the following configuration:
\begin{figure}[h]
\begin{centering}
\includegraphics{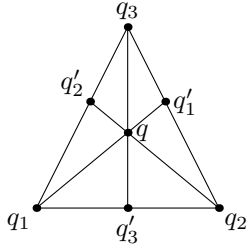}
\caption{\footnotesize{Configurations of the points $q, q_1, \ldots, q_3'$}\label{figure2}}
\end{centering}
\end{figure}

Let $\nu \colon Y \rightarrow \mathbb{P}^2$ be the blowup of these points.
Denote by $Q_i$ (respectively $Q_i'$, $Q$) the exceptional curve of $Y$ over $q_i$ (respectively $q_i'$, $q$)
and by $L$ the pullback of a general line by $\nu$.
Then
$\mathrm{Pic}(W)=\mathbb{Z}L\oplus\mathbb{Z}Q \oplus (\oplus_{i=1}^3 \mathbb{Z}Q_i \oplus \mathbb{Z}Q_i').$
The surface $Y$ has six disjoint $(-2)$-curves and  their divisor classes are:
\begin{align*}
M_i=L-Q-Q_i-Q_i',\ N_i=L-Q_i'-Q_{i+1}-Q_{i+2}\ \text{for}\ i=1, 2, 3.
\end{align*}
Let $\Lambda_i$ be the proper transform of the line $\overline{q_{i+1}'q_{i+2}'}$, i.e.,
$\Lambda_i \equiv L-Q_{i+1}'-Q_{i+2}'$ for $i=1, 2, 3$.

We describe four base-point-free pencils of rational curves on $Y$.
They are $|\Phi|:=|2L-Q-Q_1-Q_2-Q_3|$ and $|\Phi_i|:=|2L-Q-Q_i-Q_{i+1}'-Q_{i+2}'|$ ($i=1, 2, 3$).
The singular members of $|\Phi|$ are $M_1+2Q_1'+N_1, M_2+2Q_2'+N_2, M_3+2Q_3'+N_3$
and those of $|\Phi_i|$ (fixed i) are $\Lambda_i+M_i+Q_i', M_{i+1}+2Q_{i+1}+N_{i+2}, M_{i+2}+2Q_{i+2}+N_{i+1}$.
Also note that $\Phi_i+N_i \equiv -K_Y$ for $i=1, 2, 3$.

Let $\zeta \colon Y \rightarrow \Upsilon$ be the morphism
contracting the five $(-2)$-curves $M_1, M_2, N_2, M_3$ and $N_3$.
Then $\Upsilon$ has five nodes and contains a unique $(-2)$-curve $\zeta(N_1)$ in the smooth locus.

Now we define the following effective divisors on $Y$:
\begin{align}
D_1&:=\Lambda_1+\Phi_1+M_3,& D_2&:=\Lambda_2, &   D_3&:=Q_1'+\Phi+N_3, &\nonumber\\
D_{g,\i}&:=N_1+N_2,        & D_{g,-\i}&:=M_2, &   D_{g_1g, \i}&:=0, &D_{g_1g, -\i}:=M_1.  \label{eq:specialdata1}
\end{align}
We also define the following divisors:
\begin{align}
\L_{\chi}:&=4L-2Q-2Q_1-Q_2-Q_2'-Q_3-2Q_3', \label{eq:specialdata2}\\
\L_{\rho}:&=4L-2Q-2Q_1-Q_1'-2Q_2-Q_2'-Q_3-Q_3'. \nonumber
\end{align}
We require that $\Phi \in |\Phi|$ and $\Phi_1 \in |\Phi_1|$ are smooth curves
such that the divisor $D=D_1+\ldots+D_{g_1g,-\i}$ has only nodes.
These divisors satisfy \eqref{eq:Hcover} in \propref{prop:abeliancover}.
So there is a finite Galois $H$-cover $\widehat{\pi} \colon X \rightarrow Y$ and $X$ is normal.

We use \cite[Proposition~3.1 and Proposition~3.3]{pardini} to analyze the singular locus of $X$.
\begin{lem}\label{lem:singularlocus}
Let $m:=\Lambda_2 \cap M_2$ and $n:=\Lambda_2 \cap N_2$.
\begin{enumerate}[\upshape (a)]
\item The inverse image $\widehat{\pi}^{-1}(m)$ (resp. $\widehat{\pi}^{-1}(n)$)
       consists of two points $\widehat{m_1}$ and $\widehat{m_2}$ (resp. $\widehat{n_1}$ and $\widehat{n_2}$),
       each of which has stabilizer $<g>$.
\item The points $\widehat{m_1}, \widehat{m_2}, \widehat{n_1}$ and $\widehat{n_2}$ are
      exactly the singularities of $X$ and they are nodes.
\item The curve $\widehat{\pi}^{-1}(M_2)$ is a disjoint union of two irreducible smooth curves
      $\widehat{M}_{21}$ and $\widehat{M}_{22}$,
      and $\widehat{M}_{2j}$ has self intersection number $(-\frac 12)$
      and $\widehat{m}_j \in \widehat{M}_{2j}$ for $j=1, 2$.
      The curve $\widehat{\pi}^{-1}(N_2)$ consists of two irreducible smooth curves
      $\widehat{N}_{21}$ and $\widehat{N}_{22}$,
      and $\widehat{N}_{2j}$ has self intersection number $(-\frac 12)$
      and $\widehat{n}_j \in \widehat{N}_{2j}$ for $j=1, 2$.
\item The curve $\widehat{\pi}^{-1}(M_3)$ is a disjoint union four $(-1)$-curves and
      so is $\widehat{\pi}^{-1}(N_3)$.
\item The curve $\widehat{\pi}^{-1}(M_1)$ is a $(-1)$-curve.
\end{enumerate}
\end{lem}
\proof
\cite[Proposition~3.1]{pardini} shows that
$X$ is smooth outside $\widehat{\pi}^{-1}(m)$ and $\widehat{\pi}^{-1}(n)$.
Note that $M_2$ intersects only one irreducible component of $D-M_2$; that is $M_2\Lambda_2=1$.
Because $\Lambda_2 =D_2$, $M_2 \le D_{g, -\i}$ and $[H:<g>]=2$,
we conclude that $\widehat{\pi}^{-1}(m)$ consists of two points, each of which has stabilizer $<g>$.
These two points are nodes of $X$ according to \cite[Proposition~3.3]{pardini}.
For the same reason, we have $\widehat{\pi}^{-1}(M_2)=\widehat{M}_{21}\cup \widehat{M}_{22}$ with
 $\widehat{M}_{21}\cap \widehat{M}_{22}=\emptyset$ and $\widehat{\pi}|_{\widehat{M}_{2j}} \colon
\widehat{M}_{2j} \rightarrow M_2$ is an isomorphism.
We also have $\widehat{\pi}^*(M_2)=4\widehat{M}_{21}+4\widehat{M}_{22}$.
Thus (a)-(c) follow from the discussion above.
Similar arguments apply to (d) and (e).
For (d), just note that $M_3(\le D_1)$ and $N_3 (\le D_3)$ are connected irreducible components of $D$.
And (e) follows from the observation that
$M_1(=D_{g_1g, -\i})$ intersects exactly two irreducible components of $D-M_1$ and $M_1D_1=M_1D_3=1$.
\qed\smallskip

Now we explain how to obtain the smooth minimal model of $X$.
On the minimal resolution $\widetilde{X}$ of $X$,
the strict transforms of $\widehat{M}_{21}, \widehat{M}_{22}, \widehat{N}_{21}$ and $\widehat{N}_{22}$ are $(-1)$-curves.
Each of these $(-1)$-curves intersects transversely at one point with exactly one of the four $(-2)$-curves
over the nodes of $X$.
So we can contract the four curves $\widehat{M}_{21}, \widehat{M}_{22}, \widehat{N}_{21}$ and $\widehat{N}_{22}$ of $X$
to smooth points on another surface.

Let $\theta \colon X \rightarrow S$ be the morphism contracting the disjoint union of
the nine $(-1)$-curves $\widehat{\pi}^{-1}(M_3), \widehat{\pi}^{-1}(N_3), \widehat{\pi}^{-1}(M_1)$
and the four curves $\widehat{M}_{21}, \widehat{M}_{22}, \widehat{N}_{21}$ and $\widehat{N}_{22}$.
Then there is a smooth $H$-cover $\pi \colon S \rightarrow \Upsilon$
such that the outer square of the following diagram \eqref{diag:diag2} commutes.
\begin{align}\label{diag:diag2}
\xymatrix{
    X \ar_{\theta_2}"1,2" \ar@/^4mm/^{\theta}"1,3" \ar^{\widehat{\pi}_2}"2,1" \ar@/_4mm/_{\widehat{\pi}}"3,1" &
    \overline{V} \ar^{}"1,3" \ar^{\overline{\pi}}"2,2" & S \ar^{\pi}"3,3" \\
    X_1 \ar_{\delta}"2,2" \ar^{\widehat{\pi}_1}"3,1" & \overline{W} &\\
    Y \ar^{\zeta}"3,3" & &\Upsilon
}
\end{align}
We confirm that $S$ is the smooth minimal model of $X$ by the following proposition.
\begin{prop}\label{prop:specialinoue}
The surface $S$ is an Inoue surface.
\end{prop}
\proof
From \cite[Theorem~2.1]{pardini}, we obtain
$\L_{\rho^2} =2\L_\rho-D_2-D_3-D_{g, -\i}-D_{g_1g, -\i}$.
Then $2\L_{\rho^2} \equiv D_{g, \i}+D_{g, -\i}+D_{g_1g, \i}+D_{g_1g, -\i}=M_1+N_1+M_2+N_2$
by \propref{prop:abeliancover}~(b) and \eqref{eq:specialdata2}.
Let $\widehat{\pi}_1 \colon X_1 \rightarrow Y$ be the double cover
branched along the four disjoint $(-2)$-curves $M_1, N_1, M_2$ and $N_2$.
Note that $\rho^2$ is the unique character of $H^*$ which is trivial on $G$.
So the Galois group of $\widehat{\pi}_1$ is $H/G$
and the cover $\widehat{\pi}$ factors through a $G$-cover $\widehat{\pi}_2 \colon X \rightarrow X_1$.

We have $2K_{X_1}=\widehat{\pi}_1^*(2K_Y+M_1+N_1+M_2+N_2)$ and $K_{X_1}^2=0$.
The inverse images of $M_1, N_1, M_2$ and $N_2$ under $\widehat{\pi}_1$ are $(-1)$-curves.
Also $\widehat{\pi}_1^{-1}(M_3)$ is a disjoint union of two $(-2)$-curves and so is $\widehat{\pi}_1^{-1}(N_3)$.
Let $\delta \colon X_1 \rightarrow \overline{W}$ be the morphism
contracting three $(-1)$-curves $\widehat{\pi}_1^{-1}M_1, \widehat{\pi}_1^{-1}M_2$ and $\widehat{\pi}_1^{-1}N_2$.
Then $\overline{W}$ is a weak Del Pezzo surface of degree three.

Let $\theta_2 \colon X \rightarrow \overline{V}$ be the morphism contracting the curves $\widehat{\pi}^{-1}(M_1), \widehat{\pi}^{-1}(M_2)$ and $\widehat{\pi}^{-1}(N_2)$.
We obtain a smooth Galois finite $G$-cover $\overline{\pi} \colon \overline{V} \rightarrow \overline{W}$
and a commutative diagram \eqref{diag:diag2}.
The branch locus of $\overline{\pi}$ is
\begin{align}
\overline{\Delta_1}=\overline{\Lambda_1}+\overline{\Phi_1}+\overline{M_3},\
\overline{\Delta_2}=\overline{N_1}+\overline{\Lambda_2},\
\overline{\Delta_3}=\overline{Q_1'}+\overline{\Phi}+\overline{N_3} \label{eq:specialinouedata}
\end{align}
Here we denote by $\overline{\Lambda_1}=\delta\widehat{\pi}_1^{-1}(\Lambda_1)$, etc.
We claim that
\begin{enumerate}[\upshape (i)]
\item $\overline{\Lambda_1}$,  $\overline{N_1}$ and $\overline{Q_1'}$ are $(-1)$-curves;
\item $\overline{\Phi_1}$ and $\overline{\Lambda_2}$ are $0$-curves, and $\overline{\Phi}$ is a disjoint union of
      two $0$-curves in the same linear system;
\item $\overline{M_3}$ ($\overline{N_3}$) is a disjoint union of two $(-2)$-curves;
      these two $(-2)$-curves are disjoint from the $(-1)$-curves in (i);
\item $\overline{\Lambda_1} +\frac 12 \overline{\Phi}$, $\overline{N_1}+\overline{\Phi_1}$ and
      $\overline{Q_1'}+\overline{\Lambda_2}$
      and  are linearly equivalent to $-K_{\overline{W}}$.
\end{enumerate}
For example, because the general member of $|\Phi|$ is disjoint from $M_1+N_1+M_2+N_2$,
the curve $\widehat{\pi}_1^{-1}(\Phi)$ is a disjoint union of two $0$-curves in the same linear system and so is
$\overline{\Phi}$.
In particular, $K_W\overline{\Phi}=-4$ and $\frac 12 \overline{\Phi}$ is well defined in $\mathrm{Pic}(\overline{W})$.
For the $(-1)$-curve $\Lambda_1$ on $Y$,
since $\Lambda_1M_1=\Lambda_1N_1=1$ and $\Lambda_1M_2=\Lambda_1N_2=0$,
the curve $\widehat{\pi}_1^{-1}(\Lambda_1)$ is a $(-2)$-curve,
and it intersects with $\widehat{\pi}_1^{-1}(M_1)$ transversely at one point
and it is disjoint from $\widehat{\pi}_1^{-1}(M_2)$ and $\widehat{\pi}_1^{-1}(N_2)$.
So $\overline{\Lambda_1}$ is a $(-1)$-curve.
Moreover, we have
$\overline{\Lambda_1}\overline{\Phi}=\widehat{\pi}_1^*(\Lambda_1)\widehat{\pi}_1^*({\Phi})=2\Lambda_1\Phi=4$.
Finally,
the algebraic index theorem yields $\overline{\Lambda_1}+\frac  12 \overline{\Phi} \equiv -K_W$.
Other statements can be proved in the same manner.

Comparing \eqref{eq:specialinouedata} to \eqref{eq:inouedata},
we conclude that $S$ is an Inoue surface.
\qed.\smallskip

When $\Phi$ and $\Phi_1$ vary, we obtain a $2$-dimensional family of Inoue surfaces with automorphism groups
isomorphic to $\mathbb{Z}_2 \times \mathbb{Z}_4$.

\begin{rem}
We may directly show that $K_S$ is ample, $K_S^2=7$ and $p_g(S)=0$ for the surface $S$ in \eqref{diag:diag2}.
According to the proof of \cite[Proposition~4.2]{pardini},
we have
\begin{align*}
4K_X&=\widehat{\pi}^*(4K_Y+2D_1+2D_2+2D_3+3D_{g,\i}+3D_{g,-\i}+3D_{g_1g, \i}+3D_{g_1g, -\i})\\
    &=\widehat{\pi}^*(-K_Y+\Phi+\Phi_1)+\widehat{\pi}^*(M_1+2M_2+2N_2+2M_3+2N_3)
\end{align*}
It follows that $4K_S=\pi^*(-K_\Upsilon+\phi+\phi_1)$ and $K_S^2=7$,
where $|\phi|$ and $|\phi_1|$ are base-point-free pencils on $\Upsilon$ induced by $|\Phi|$ and $|\Phi_1|$.
The linear system $|-K_Y+\Phi+\Phi_1|$ is base point free,
and the corresponding morphism contracts exactly the nodal curves $M_1, M_2, N_2, M_3$ and $N_3$.
Hence $|-K_\Upsilon+\phi+\phi_1|$ is ample and so is $K_S$.
For each $\psi \in H^*$, we can calculate $\L_\psi$ by \cite[Theorem~2.1]{pardini}
and then easily show that $H^0(Y, \O_Y(K_Y+\L_\psi))=0$.
It follows that $p_g(S)=p_g(X)=0$ by \cite[Proposition~4.1]{pardini}.
\end{rem}

\begin{rem}We remark that \thmref{thm:inoueauto} contributes to the study of the moduli space of the Inoue surfaces.
Let $\mathcal{M}_{1,7}^{\mathrm{can}}$ be the Gieseker moduli space of canonical models of surfaces of general type
with $\chi(\O)=1$ and $K^2=7$ (cf.~\cite{gieseker}).
Let $S$ be any Inoue surface.
Denote by $[S]$ the corresponding point in $\mathcal{M}_{1,7}^{\mathrm{can}}$ and
by $\mathrm{B}(S)$ be the base of the Kuranishi family of deformations of $S$.
Recall the facts that  the tangent space of $\mathrm{B}(S)$ is $H^1(S, \Theta_S)$,
where $\Theta_S$ is the tangent sheaf of $S$,
and that the germ $(\mathcal{M}_{1,7}^{\mathrm{can}}, [S])$ is analytically isomorphic to
$\mathrm{B}(S)/\mathrm{Aut}(S)$.
It has been shown in \cite{inouemfd} that the group $G$ acts trivially on $H^1(S, \Theta_S)$
and $\mathrm{B}(S)$ is smooth of dimension $4$.

Now assume that $S$ is a special Inoue surface constructed here.
We can use the same method as in the proof of \cite[Theorem~5.1]{inouemfd} to conclude that
the invariant subspace of $H^1(S, \Theta_S)$ for the $H$-action has dimension $2$.
Note that $\mathrm{Aut}(S)=H$ and $H/G \cong \mathbb{Z}_2$.
Combining the result of \cite{inouemfd},
we see that $(\mathcal{M}_{1,7}^{\mathrm{can}}, [S])$ is analytically isomorphic to
$(\mathbb{C}^2 \times \mathrm{Spec}\ \mathbb{C}[x,y,z]/(xz-y^2), 0)$.

\end{rem}

\paragraph{Acknowledgement.}
The author thanks Xiaotao~Sun for all the support during the preparation of the article.
The author is greatly indebted to Jinxing~Cai for many discussions and helpful suggestions.
This work is partially supported by the China Postdoctoral Science Foundation (Grant No.: 2013M541062).

\end{document}